\documentclass[12pt,fceqn]{article}
\usepackage{amsmath}
\usepackage{amsfonts}
\usepackage{cite}
\usepackage{amssymb,graphics,graphicx,subfigure,caption2,rotating}
\usepackage{amsmath, latexsym}
\usepackage[amsmath,thmmarks]{ntheorem}
\numberwithin{equation}{section}
\newtheorem{theorem}{Theorem}[section]
\newtheorem{lemma}{Lemma}[section]

\newtheorem{prop}{Proposition}[section]

\newcommand{\beq}{\begin{equation}}
\newcommand{\eeq}{\end{equation}}
\newcommand{\beqn}{\begin{eqnarray}}
\newcommand{\eeqn}{\end{eqnarray}}
\date{}
\begin{document}

\date{}
\title{New error bounds for the extended vertical LCP\thanks{This research was
supported by National Natural Science Foundation of China
(No.11961082).}}
\author{Shiliang Wu\thanks{Corresponding author: slwuynnu@126.com}, Hehui Wang\thanks{wanghehui1994@126.com}\\
{\small{\it $^{\dag}$School of Mathematics, Yunnan Normal University,}}\\
{\small{\it Kunming, Yunnan, 650500, P.R. China}}\\
{\small{\it $^{\ddag}$School of Mathematics and Statistics, Yunnan University,}}\\
{\small{\it Kunming, Yunnan, 650091, P.R. China}}\\
}
 \maketitle
\begin{abstract}
In this paper, by making use of this fact that for $a_{j}, b_{j}\in
\mathbb{R}$, $j=1,2,\ldots,n$, there are $\lambda_{j}\in [0,1]$ with
$\sum_{j=1}^{n}\lambda_{j}=1$ such that
\[
\min_{1\leq j\leq n}\{a_{j}\}-\min_{1\leq j\leq
n}\{b_{j}\}=\sum_{j=1}^{n}\lambda_{j}(a_{j}-b_{j}),
\]
 some new error bounds of the extended vertical LCP under the row
$\mathcal{W}$-property are obtained, which cover the error bounds in
[Math. Program., 106 (2006) 513-525] and [Comput. Optim. Appl., 42
(2009) 335-352]. Not only that, these new error bounds skillfully
avoid the inconvenience caused by the row rearrangement technique
for error bounds to achieve the goal of reducing the computation
workload, which was introduced in the latter paper mentioned above.
Besides, with respect to the row $\mathcal{W}$-property, two new
sufficient and necessary conditions are obtained.


\end{abstract}

\section{Introduction}
For $A_{j}\in \mathbb{R}^{n\times n}$ and $q_{j}\in \mathbb{R}^{n}$
($j=0,1,2,\ldots,k$) being given known matrices and the source
terms, the extended vertical linear complementarity problem is to
find $x\in \mathbb{R}^{n}$ such that
\begin{equation}\label{eq:11}
r(x):=\min\{A_{0}x+q_{0}, A_{1}x+q_{1},\ldots, A_{k}x+q_{k}\}=0,
\end{equation}
where  $\min$ is the component minimum operator. Here, Eq.
(\ref{eq:11}) is denoted by EVLCP($\mathbf{A},\mathbf{q}$) for
short, where
\[
\mathbf{A}=(A_{0},A_{1},\ldots, A_{k})\ \mbox{and} \
\mathbf{q}=(q_{0},q_{1},\ldots, q_{k}).
\]
When $A_{0}=I$ and $q_{0}=0$ in (\ref{eq:11}), where $I$ denotes the
identity matrix, the EVLCP ($\mathbf{A},\mathbf{q}$) reduces to the
vertical LCP, which was introduced by Cottle and Dantzig
\cite{Cottle70}, also see \cite{Gowda}. Further, when $A_{0}=I$ and
$q_{0}=0$ and $k=1$ in (\ref{eq:11}), the EVLCP
($\mathbf{A},\mathbf{q}$) comes back to the standard  LCP
($A_{1},q_{1}$), see \cite{Cottle92, Murty88}.

So far,  it has been found that the EVLCP ($\mathbf{A},\mathbf{q}$)
has more and more applications in the many fields, like, such as
nonlinear networks \cite{Fujisawa72}, control theory \cite{Sun89},
the mixed lubrication problem \cite{Oh86}, stochastic impulse
control games \cite{Zabaljauregui21},  the boundary value problem
\cite{Goeleven96}, the generalized bimatrix games \cite{Gowda96},
the generalized Leontief input-output model \cite{Ebiefung93}, the
discrete HJB equations \cite{Zhou08}, volterra ecosystem
\cite{Habetler92} and so on. There exist many literatures to pay
attention to the existence of solutions and algorithms for the
EVLCP ($\mathbf{A},\mathbf{q}$), see \cite{Gowda, Qi99, Sun89,
Gowda962, Mezzadri,Zabaljauregui21}.

For the EVLCP ($\mathbf{A},\mathbf{q}$), another important and
interested topic in theory is error bound, which has drawn
widespread attention because of an important tool in theoretical
analysis, including convergence analysis, sensitive analysis and
verification of the solutions. At present, there have some results
about error bounds in the literature. For instance, assume that
$x^{\ast}$ is the unique solution of the LCP ($A_{1},q_{1}$), for
$A_{1}$ being $P$-matrix, a well-know global error bound was given
in \cite{Mathias90} by Mathias and Pang, and described below
\begin{equation}\label{eq:12}
\|x-x^{\ast}\|_{\infty}\leq\frac{1+\|A_{1}\|_{\infty}}{\alpha(A_{1})}\|r(x)\|_{\infty},
\ \mbox{for any }\ x\in \mathbb{R}^{n},
\end{equation}
where
\[
\alpha(A_{1}):=\min_{\|x\|_{\infty}=1}\big\{\max_{1\leq i\leq
n}x_{i}(A_{1}x)_{i}\big\}.
\]
By the equivalent form of the minimum function, Chen and Xiang in
\cite{Chen06} obtained the following error bound in $\|\cdot\|_{p}\
(p\geq1, \mbox{or} \ p=\infty)$ norms,
\begin{equation}\label{eq:13}
\|x-x^{\ast}\|_{p}\leq\max_{d\in
[0,1]^{n}}\|(I-D+DA_{1})^{-1}\|_{p}\|r(x)\|_{p}, \ \mbox{for  any }
\ x\in \mathbb{R}^{n},
\end{equation}
where $D=\mbox{diag}(d)$ with $d\in [0,1]^{n}$, which is sharper
than (\ref{eq:12}) in $\|\cdot\|_{\infty}$, see \cite{Chen06}.
Moreover, for $A_{1}$ being an $H_{+}$-matrix, Chen and Xiang
confirmed
\begin{equation}\label{eq:14}
\max_{d\in [0,1]^{n}}\|(I-D+DA_{1})^{-1}\|_{p}\leq \|\langle
A_{1}\rangle ^{-1}\max(\wedge_{1},I)\rangle\|_{p},
\end{equation}
where $\wedge_{1}$ is the diagonal part of $A_{1}$ and $\langle
A_{1}\rangle$ is its comparison matrix (i.e., $\langle
A_{1}\rangle_{ii}=|(A_{1})_{ii}|, \langle
A_{1}\rangle_{ij}=-|(A_{1})_{ij}|$ for $i\neq j$).

Recently, by the row rearrangement technique introduced by Zhang et
al. in \cite{Zhang09}, i.e., $\textbf{A}'=(A'_{0}, A'_{1}, \ldots,
A'_{k})$ is called a row rearrangement of
$\mathbf{A}=(A_{0},A_{1},\ldots, A_{k})$ if for any $i\in
N:=\{1,2,\ldots,n\}$,
\[
(\textbf{A}'_{j})_{i\cdot}=(A_{ji})_{i\cdot}\in
\{(A_{0})_{i\cdot},(A_{1})_{i\cdot},\ldots,
(A_{k})_{i\cdot}\}=\{(A'_{0})_{i\cdot},(A'_{1})_{i\cdot},\ldots,
(A'_{k})_{i\cdot}\},
\]
where $(\cdot)_{i\cdot}$ means the i-th row of a given matrix, for
the EVLCP ($\mathbf{A},\mathbf{q}$), assume that
$\mathbf{A}=(A_{0},A_{1},\ldots, A_{k})$ has the row
$\mathcal{W}$-property (see the following definition), Zhang et al.
in \cite{Zhang09} presented the following result,
\begin{equation}\label{eq:15}
\|x-x^{\ast}\|\leq\alpha(\textbf{A})\|r(x)\|, \ \mbox{for any } \
x\in \mathbb{R}^{n},
\end{equation}
where
\begin{equation}\label{eq:16}
\alpha(\textbf{A})=\max_{\textbf{A}'\in
\mathcal{R}(\mathbf{A})}\max_{j<l\in\{0,1,\ldots,k\}}\max_{d\in
[0,1]^{n}}\|((I-D)A_{j}'+A'_{l})^{-1}\|,
\end{equation}
$A_{j}', A'_{l}\in \mathbb{R}^{n\times n}$ are any two blocks in
$\textbf{A}'\in \mathcal{R}(\mathbf{A})$ with
$\mathcal{R}(\mathbf{A})$ being the set of all row rearrangements of
$\mathbf{A}$. When the norm in (\ref{eq:15}) is taken as
$\|\cdot\|_{\infty}$, we denote
\[
\alpha_{\infty}(\textbf{A})=\max_{\textbf{A}'\in
\mathcal{R}(\mathbf{A})}\max_{j<l\in\{0,1,\ldots,k\}}\max_{d\in
[0,1]^{n}}\|((I-D)A_{j}'+A'_{l})^{-1}\|_{\infty}.
\]

There is no doubt that  Eq. (\ref{eq:15}) presents a general result
for the upper global  error bound of the EVLCP
($\mathbf{A},\mathbf{q}$) under the row $\mathcal{W}$-property. By
investigating Eq. (\ref{eq:15}), clearly, it is easy to know that
before obtaining Eq. (\ref{eq:15}), we have to calculate and obtain
the exact value of (\ref{eq:16}). Whereas, in the implementations,
calculating the exact value of (\ref{eq:16}) is a very difficult
task because there is a rearrangement for the row of matrix
$\mathbf{A}$, in particular, when the order of matrix $\mathbf{A}$
is large. It is a main motivation of this present paper. In this
paper, to overcome this disadvantage caused by the row rearrangement
technique in essence, we have to carve out a new approach to obtain
the error bound for the EVLCP ($\mathbf{A},\mathbf{q}$). Our
approach inspired by the work in \cite{Chen06}, we first develop a
general equivalent form of the minimum function. Then, based on
this, some new error bounds for the EVLCP ($\mathbf{A},\mathbf{q}$)
are obtained. Not only that, these new error bounds cover some
existing results in \cite{Chen06} and \cite{Zhang09} as well.
Meanwhile, it avoids the row rearrangement of the system matrix
$\mathbf{A}$. Incidentally, for the row $\mathcal{W}$-property,  two
new sufficient and necessary conditions are given.

The rest of the article expands as follows. First, from the view of
the calculation time, section 2 further discusses the result in
(\ref{eq:15}). Secondly, section 3 provides some error bounds of the
EVLCP ($\mathbf{A},\mathbf{q}$) under the row $\mathcal{W}$-property
by a general equivalent form of the minimum function, which is
entirely different from Eq. (\ref{eq:15}). Thirdly, section 4
presents some numerical examples to show the feasibility of the
error bound. Finally, section 5 is a brief conclusion.

By the way, the following notations,  definitions and results will
be used throughout the paper, which can be founded in \cite{Berman,
Gowda}. Let $A = (a_{ij})$ and $B = (b_{ij}) \in \mathbb{R}^{n\times
n}$. Then $A\geq (>)B$ means $a_{ij}\geq (>) b_{ij}$ for
$i,j=1,2,\ldots,n$. We indicate $|A|=(|a_{ij}|)$. Matrix $A =
(a_{ij})$ is called a strictly diagonal dominant  matrix if
$|a_{ii}|>\sum_{j\neq i}|a_{ij}|,\ i\in N:=\{1,2,\ldots,n\}$.
$\rho(\cdot)$ indicates the spectral radius of the matrix. A block
matrix $\mathbf{A}=(A_{0},A_{1},\ldots, A_{k})$ has the row
$\mathcal{W}$-property if
\[
\min(A_{0}x, A_{1}x,\ldots, A_{k}x)\leq0\leq \max(A_{0}x,
A_{1}x,\ldots, A_{k}x)\Rightarrow x=0.
\]
The EVLCP ($\mathbf{A},\mathbf{q}$) has a unique solution for any
$\mathbf{q}$ if and only if $\mathbf{A}$ has the row
$\mathcal{W}$-property.

\section{Old error bound}
In this section, we will discuss the result in (\ref{eq:15}) from
the angle of the calculation time.

To calculate $\alpha(\textbf{A})$, a natural question is how many
times you need to calculate $\|((I-D)A_{j}'+A'_{l})^{-1}\|$.
Further, when confronting Eq. (\ref{eq:16}), the first thought is to
calculate the number of elements in $\mathcal{R}(\mathbf{A})$, that
is, we need to know that how many the row rearrangements of
$\mathbf{A}=(A_{0},A_{1},\ldots, A_{k})$ there are in
$\mathcal{R}(\mathbf{A})$. To answer this question, we get
Proposition 2.1.

\begin{prop}
Let $\mathbf{A}=(A_{0},A_{1},\ldots, A_{k})$ with $A_{j}\in
\mathbb{R}^{n\times n}$ $(j=0,1,\ldots,k)$. Then the cardinality of
$\mathcal{R}(\mathbf{A})$ is $[(k+1)!]^{n}$.
\end{prop}
\textbf{Proof.}  If only one row is rearranged, noting that a total
number of matrices are $k+1$, there are a total of $(k+1)!$
different sorting methods. Because there are a total of $n$ rows and
the arrangement of different rows can be arbitrarily combined, there
are a total of $[(k+1)!]^{n}$ different sorting methods. Therefore,
matrix $\mathbf{A}$ has $[(k+1)!]^{n}$ different row rearrangements,
that is, the cardinality of $\mathcal{R}(\mathbf{A})$ is $[(
k+1)!]^{n}$. $\hfill{} \Box$

In fact, once a certain $\textbf{A}'$ is selected, there exist the
total $C_{k+1}^{2}=\frac{k(k+1 )}{2}$ kinds of options for choosing
$A_{j}', A'_{l}$. From Proposition 2.1, a total of elements in
$\mathcal{R}(\mathbf{A})$ are $[(k+1)!]^{n}$, so the calculation
times of $\|((I-D)A_{j}'+A'_{ l})^{-1}\|$ are at most
$\frac{k(k+1)}{2}[(k+1)!]^{n}$ times to obtain $\alpha(\textbf{A})$.
Therefore, the calculation times of Theorem 4.1 in \cite{Zhang09}
for $k=1$ achieve at most $2^{n}$ times. Whereas,  the calculation
time of  Theorem 3.1 in \cite{Zhang09} is 1. The next question is
whether Theorem 4.1 in \cite{Zhang09} for $k=1$ is in contradiction
with Theorem 3.1 \cite{Zhang09} because the form of Theorem 4.1 in
\cite{Zhang09} for $k=1$ is different from in Theorem 3.1 in
\cite{Zhang09} as well. To answer this question, we require
Proposition 2.2.

\begin{prop}
Let $\mathbf{A}=(A_{0},A_{1})$ with $A_{0}, A_{1}\in
\mathbb{R}^{n\times n}$. Then
\[
\max_{\textbf{A}'\in \mathcal{R}(\mathbf{A})}\max_{d\in
[0,1]^{n}}\|((I-D)A'_{0}+A'_{1})^{-1}\|= \max_{d\in
[0,1]^{n}}\|((I-D)A_{0}+DA_{1})^{-1}\| .
\]
\end{prop}
\textbf{Proof.} For any $\textbf{A}'=(A_{0}', A'_{1})$,
$D=\mbox{diag}(d)$, $d\in [0,1]^{n}$, $i=1,2,\cdots,n$, we have
\[
((I-D)A'_{0}+DA'_{1})_{i.}=(1-d_i)(A'_{0})_{i.}+d_i(A'_{1 })_{i.}.
\]
According to the definition of row rearrangement, we get
\[
\{(A'_{0})_{i.},(A'_{1})_{i.}\}=\{(A_{0})_{i.},(A_{1})_{i.}\}.
\]
When $(A'_{0})_{i.}=(A_{0})_{i.}$, $(A'_{1})_{i.}=(A_{1})_{i.}$, we
take $\hat{d}_i=d_i$; when $(A'_{0})_{i.}=(A_{1})_{i.}$,
$(A'_{1})_{i.}=(A_{0})_{i.}$, we take $\hat{d}_i=1-d_i$,
$\hat{d}_i\in [0,1]$. Further, we have
\[
(1-d_i)(A'_{0})_{i.}+d_i(A'_{1})_{i.}=(1-\hat{d}_i)(A_{0})_{i.}+\hat{d}_i(A_{1})_{i.}.
\]
Let $\hat{D}=\mbox{diag}(\hat{d})$. Then
\[
(I-D)A'_{0}+DA'_{1}=(I-\hat{D})A_{0}+\hat{D}A_{1}.
\]
Therefore, from the arbitrariness of $\textbf{A}'$, we obtain
\[
\max_{\textbf{A}'\in \mathcal{R}(\mathbf{A})}\max_{d\in
[0,1]^{n}}\|((I-D)A'_{0}+A'_{1})^{-1}\|= \max_{d\in
[0,1]^{n}}\|((I-D)A_{0}+DA_{1})^{-1}\|.
\]
This completes the proof. $\hfill{} \Box$

Based on Proposition 2.2, in essence,  Theorem 4.1 in \cite{Zhang09}
for $k=1$ is in line with Theorem 3.1 \cite{Zhang09}.  The
calculation times of both are the same, 1 time.

According to the process of the proof of Proposition 2.2, we can see
that when calculating $\|((I-D)A_{j}'+DA'_{l})^{-1}\|$, as long as
two of the $k+1$ elements in each row are chosen, no matter how many
$A_{j}, A'_{l}$ they can be combined into, ultimately, their
calculation results are the same, so there is no need to reconsider.
Therefore, calculating $\alpha(\textbf{A})$ only needs to calculate
$(C_{k+1}^2)^{n}$ times $\|((I-D)A'_{j}+A'_{l})^{-1}\|$. From this,
we present Proposition 2.3.

\begin{prop}
Let $\mathbf{A}=(A_{0},A_{1},\ldots, A_{k})$. Then there needs to
calculate $\alpha(\emph{\textbf{A}})$  at most $(C_{k+1 }^2)^n$
times $\|((I-D)A_{j}'+DA'_{l})^{-1}\|$, i.e.,
\begin{equation}
\alpha(\emph{\textbf{A}})=\max_{\textbf{A}'\in
\mathcal{R}(\mathbf{A})}\max_{B_{1},B_{2}}\max_{d\in
[0,1]^{n}}\|((I-D)B_{1}+DB_{2})^{-1}\|,
\end{equation} where
$(B_{1})_{i.}=(A_{j})_{i.}$, $(B_{2})_{i.}=(A_{l})_{i. }$,
$j<l\in\{0,1,\ldots,k\}$.
\end{prop}

To further explain Proposition 2.3, Example 2.1 is provided.

\textbf{Example 2.1}  Let $\mathbf{A}=(A_{0}, A_{1}, A_{2})$, where
\[
A_{0}=\left[\begin{array}{cc}
1&1\\
-1&1
\end{array}\right],\ A_{1}=\left[\begin{array}{cc}
1&0\\
-2&1
\end{array}\right], \ A_{2}=\left[\begin{array}{cc}
2&1\\
0&1
\end{array}\right].
\]

When using Eq. (2.1), from Proposition 2.3, its calculation times
are \textbf{9}. Specifically as follows:

Let
\[B_{1}=\left[\begin{array}{cc}
1&1\\
-1&1
\end{array}\right],
 B_{2}=\left[\begin{array}{cc}
1&0\\
-2&1
\end{array}\right].
\]
Then
\begin{align*}
\mu_{1}=&\max_{d\in
[0,1]^{n}}\|((I-D)B_{1}+DB_{2})^{-1}\|_{\infty}=3;
\end{align*}
let
\[B_{1}=\left[\begin{array}{cc}
1&1\\
-1&1
\end{array}\right],
 B_{2}=\left[\begin{array}{cc}
2&1\\
0&1
\end{array}\right].
\]
Then
\begin{align*}
\mu_{2}=&\max_{d\in
[0,1]^{n}}\|((I-D)B_{1}+DB_{2})^{-1}\|_{\infty}=2;
\end{align*}
let
\[B_{1}=\left[\begin{array}{cc}
1&0\\
-2&1
\end{array}\right],
B_{2}=\left[\begin{array}{cc}
2&1\\
0&1
\end{array}\right].
\]
Then
\begin{align*}
\mu_{3}=&\max_{d\in
[0,1]^{n}}\|((I-D)B_{1}+DB_{2})^{-1}\|_{\infty}=3;
\end{align*}
let
\[B_{1}=\left[\begin{array}{cc}
1&1\\
-1&1
\end{array}\right],
B_{2}=\left[\begin{array}{cc}
1&0\\
0&1
\end{array}\right].
\] Then
\begin{align*}
\mu_{4}=&\max_{d\in
[0,1]^{n}}\|((I-D)B_{1}+DB_{2})^{-1}\|_{\infty}=2;
\end{align*}
let
\[B_{1}=\left[\begin{array}{cc}
1&1\\
-2&1
\end{array}\right],
B_{2}=\left[\begin{array}{cc}
1&0\\
0&1
\end{array}\right].
\] Then
\begin{align*}
\mu_{5}=&\max_{d\in
[0,1]^{n}}\|((I-D)B_{1}+DB_{2})^{-1}\|_{\infty}=3;
\end{align*}
let
\[B_{1}=\left[\begin{array}{cc}
1&1\\
-1&1
\end{array}\right],
B_{2}=\left[\begin{array}{cc}
2&1\\
-2&1
\end{array}\right].
\] Then
\begin{align*}
\mu_{6}=&\max_{d\in
[0,1]^{n}}\|((I-D)B_{1}+DB_{2})^{-1}\|_{\infty}=\frac{3}{2};
\end{align*}
let
\[B_{1}=\left[\begin{array}{cc}
1&1\\
-2&1
\end{array}\right],
 B_{2}=\left[\begin{array}{cc}
2&1\\
0&1
\end{array}\right].
\] Then
\begin{align*}
\mu_{7}=&\max_{d\in
[0,1]^{n}}\|((I-D)B_{1}+DB_{2})^{-1}\|_{\infty}=2;
\end{align*}
let
\[B_{1}=\left[\begin{array}{cc}
1&0\\
-1&1
\end{array}\right],
 B_{2}=\left[\begin{array}{cc}
2&1\\
-2&1
\end{array}\right].\]
 Then
\begin{align*}
\mu_{8}=&\max_{d\in
[0,1]^{n}}\|((I-D)B_{1}+DB_{2})^{-1}\|_{\infty}=3;
\end{align*}
let
\[B_{1}=\left[\begin{array}{cc}
1&0\\
-1&1
\end{array}\right],
 B_{2}=\left[\begin{array}{cc}
2&1\\
0&1
\end{array}\right].\]
Then
\begin{align*}
\mu_{9}=&\max_{d\in
[0,1]^{n}}\|((I-D)B_{1}+DB_{2})^{-1}\|_{\infty}=2.
\end{align*}
Based on the above computational results, we obtain
\[
\alpha_{\infty}(\textbf{A})=\max\{\mu_{1},\mu_{2},\cdots,\mu_{9}\}=3.
\]

However, directly  using Eq. (\ref{eq:15}), its calculation times
are \textbf{$\frac{k(k+1)}{2}[(k+1)!]^{n}=108$}, this is because a
lot of repetitive work are done. For example, let,
\[A'_{0}=\left[\begin{array}{cc}
1&1\\
-2&1
\end{array}\right],
 A'_{1}=\left[\begin{array}{cc}
1&0\\
-1&1
\end{array}\right],
 A'_{2}=\left[\begin{array}{cc}
2&1\\
0&1
\end{array}\right],
\]
then $A'=(A'_0,A'_1,A'_2)\in \mathcal{R}(\mathbf{A})$. Similarly,
let,
\[A''_{0}=\left[\begin{array}{cc}
1&0\\
-2&1
\end{array}\right],
 A''_{1}=\left[\begin{array}{cc}
1&1\\
-1&1
\end{array}\right],
 A''_{2}=\left[\begin{array}{cc}
2&1\\
0&1
\end{array}\right],
\]
then $A''=(A''_0,A''_1,A''_2)\in \mathcal{R}(\mathbf{A})$. Of
course, $A=(A_0,A_1,A_2)\in \mathcal{R}(\mathbf{A})$. When using Eq.
(\ref{eq:15}), we have to calculate
\[
\max_{d\in [0,1]^{n}}\|((I-D)A_{0}+DA_{1})^{-1}\|_{\infty},
\max_{d\in [0,1]^{n}}\|((I-D)A'_{0}+DA'_{1})^{-1}\|_{\infty} \]
\mbox{and}
\[\max_{d\in
[0,1]^{n}}\|((I-D)A''_{0}+DA''_{1})^{-1}\|_{\infty}.
\]
Whereas, from the proof of Proposition 2.2, the above three formulas
are equivalent, so they belong to the repeated calculation.
Proposition 2.3 just avoids this.

Example 2.1 is  simple, but it tells us that for any block matrix
$\mathbf{A}$ containing three two-by-two matrices, i.e., $k=2$ and
$n=2$, once the row rearrangement of $\mathbf{A}$ is required,  its
calculation times are at most \textbf{9} to obtain the corresponding
error bound. It's easy to imagine that for the block matrix
$\mathbf{A}$ containing $k+1$ $n$-by-$n$ matrices, making use of
Proposition 2.3 to compute the error bound of the EVLCP
($\mathbf{A},\mathbf{q}$), once the row rearrangement occurs, the
corresponding computational cost should be highly expensive and
unacceptable for the sufficiently large $k$ or  $n$! Therefore, to
face this headwind, we have to exploit a new and effective tool to
obtain the error bound of the EVLCP ($\mathbf{A},\mathbf{q}$).

\section{New error bound}

In this section, to avoid Proposition 2.3, we will give some new
error bounds for the EVLCP ($\mathbf{A},\mathbf{q}$). For this goal,
we first present some requisite lemmas.

\begin{lemma}
Let all $a_{j}, b_{j}\in \mathbb{R}$, $j=1,2,\ldots,n$. Then there
are $\lambda_{j}$ with $\lambda_{j}\in [0,1]$ such that
\begin{equation}\label{eq:31}
\min_{1\leq j\leq n}\{a_{j}\}-\min_{1\leq j\leq
n}\{b_{j}\}=\sum_{j=1}^{n}\lambda_{j}(a_{j}-b_{j}).
\end{equation}
\end{lemma}
\textbf{Proof.}  The result in Lemma 3.1 is given directly from the
mean value theorem of Lipschitz functions with the generalized
gradient. $\hfill{} \Box$

\begin{lemma}
Let $a_{j}, b_{j}, t_{j}\in \mathbb{R}$ with $a_{j}>0$, $t_{j}\in
[0,1]$, $(j=1,2,\ldots,n)$ and $\sum_{j=1}^{n}t_{j}=1$. Then
\begin{equation}\label{eq:33}
\frac{\sum_{j=1}^{n}t_{j}b_{j}}{\sum_{j=1}^{n}t_{j}a_{j}}\leq\max_{1\leq
j\leq n}\bigg\{\frac{|b_{j}|}{a_{j}}\bigg\}.
\end{equation}
In addition, if $b_{j}=1$, then the inequality (\ref{eq:33})
simplifies as
\begin{equation*}
\frac{1}{\sum_{j=1}^{n}t_{j}a_{j}}\leq\max_{1\leq j\leq
n}\bigg\{\frac{1}{a_{j}}\bigg\}.
\end{equation*}
\end{lemma}
\textbf{Proof.} The proof is straightforward. $\hfill{} \Box$

In addition, it is noted that the following inequality is still
true, i.e.,
\begin{equation*}
\sum_{j=2}^{n}t_{j}b_{j}\leq\max_{2\leq j\leq n}\{|b_{j}|\},
\end{equation*}
where  $t_{j}\in [0,1]$ and $\sum_{j=2}^{n}t_{j}\leq1$. In fact, by
the simple calculation, we have
\begin{align*}
\sum_{j=2}^{n}t_{j}b_{j}=&\frac{\sum_{j=2}^{n}t_{j}b_{j}}{1-\sum_{j=2}^{n}t_{j}+\sum_{j=2}^{n}t_{j}}\\
=&\frac{\sum_{j=2}^{n}t_{j}b_{j}}{\frac{1-\sum_{j=2}^{n}t_{j}}{n-2}(n-2)+\sum_{j=2}^{n}t_{j}}\\
=&\frac{t_{2}b_{2}+t_{3}b_{3}+\ldots+t_{n}b_{n}}{(\frac{1-\sum_{j=2}^{n}t_{j}}{n-2}+t_{2})+(\frac{1-\sum_{j=2}^{n}t_{j}}{n-2}+t_{3})+\ldots+(\frac{1-\sum_{j=2}^{n}t_{j}}{n-2}+t_{n})}\\
\leq&\frac{t_{2}|b_{2}|+t_{3}|b_{3}|+\ldots+t_{n}|b_{n}|}{(\frac{1-\sum_{j=2}^{n}t_{j}}{n-2}+t_{2})+(\frac{1-\sum_{j=2}^{n}t_{j}}{n-2}+t_{3})+\ldots+(\frac{1-\sum_{j=2}^{n}t_{j}}{n-2}+t_{n})}\\
\leq&\frac{(\frac{1-\sum_{j=2}^{n}t_{j}}{n-2}+t_{2})|b_{2}|+(\frac{1-\sum_{j=2}^{n}t_{j}}{n-2}+t_{3})|b_{3}|+\ldots+(\frac{1-\sum_{j=2}^{n}t_{j}}{n-2}+t_{n})|b_{n}|}
{(\frac{1-\sum_{j=2}^{n}t_{j}}{n-2}+t_{2})+(\frac{1-\sum_{j=2}^{n}t_{j}}{n-2}+t_{3})+\ldots+(\frac{1-\sum_{j=2}^{n}t_{j}}{n-2}+t_{n})}\\
\leq&\max_{2\leq j\leq n}\{|b_{j}|\}.
\end{align*}

\begin{lemma} \emph{\cite{Sznajder95}}
Matrix $\mathbf{A}=(A_{0},A_{1},\ldots, A_{k})$ has the row
$\mathcal{W}$-property if and only if for arbitrary nonnegative
 diagonal matrix $X_{0}, X_{1},\ldots, X_{k}$ with
$\mbox{diag}(X_{0}+X_{1}+\ldots+X_{k})>0$,
\[
\det(X_{0}A_{0}+X_{1}A_{1}+\ldots+X_{k}A_{k})\neq0.
\]
\end{lemma}

To obtain the error bound for the EVLCP ($\mathbf{A},\mathbf{q}$),
we require a new  sufficient and necessary for the  row
$\mathcal{W}$-property, see Lemma 3.4.

\begin{lemma}
Matrix $\mathbf{A}=(A_{0},A_{1},\ldots, A_{k})$ has the row
$\mathcal{W}$-property if and only if
$D_{0}A_{0}+D_{1}A_{1}+\ldots+D_{k}A_{k}$ is nonsingular for
arbitrary nonnegative diagonal matrices $D_{j}=\mbox{diag}(d_{j})$
with  $d_{j}\in [0,1]^{n}$, $(j=0,1,\ldots,k)$ and
$\sum_{j=0}^{k}D_{j}=I$.
\end{lemma}

\textbf{Proof.} The proof is straightforward by making use of Lemma
3.3. $\hfill{} \Box$

%By Lemma 3.3, matrix $\mathbf{A}=(A_{0},A_{1},\ldots, A_{k})$ has
%the row $\mathcal{W}$-property if and only if
%$X_{0}A_{0}+X_{1}A_{1}+\ldots+X_{k}A_{k}$ is nonsingular for
%arbitrary nonnegative diagonal matrix $X_{0}, X_{1},\ldots, X_{k}$
%with $\mbox{diag}(X_{0}+X_{1}+\ldots+X_{k})>0$. Noting that
%$\mbox{diag}(X_{0}+X_{1}+\ldots+X_{k})>0$, let $X=X_{0}+
%X_{1}+\ldots+X_{k}$. Then matrix $\mathbf{A}=(A_{0}, A_{1},\ldots,
%A_{k})$ has the row $\mathcal{W}$-property if and only if
%$X^{-1}(X_{0}A_{0}+X_{1}A_{1}+\ldots+X_{k}A_{k})$ is nonsingular for
%arbitrary nonnegative diagonal matrix $X_{0}, X_{1},\ldots, X_{k}$
%with $\mbox{diag}(X)>0$. Denote
%\[
%D_{0}=X^{-1}X_{0}, D_{1}=X^{-1}X_{1},\ldots, D_{k}=X^{-1}X_{k}.
%\]
%Clearly, $\mathbf{A}=(A_{0},A_{1},\ldots, A_{k})$ has the row
%$\mathcal{W}$-property if and only if
%$D_{0}A_{0}+D_{1}A_{1}+\ldots+D_{k}A_{k}$ is nonsingular for
%arbitrary nonnegative diagonal matrices $D_{j}=\mbox{diag}(d_{j})$
%with  $d_{j}\in [0,1]^{n}$, $(j=0,1,\ldots,k)$ and
%$\sum_{j=0}^{k}D_{j}=I$. $\hfill{} \Box$

\textbf{Remark 3.1} When $k=1$ in Lemma 3.4, the result in Lemma 3.4
goes back to Lemma 2.2 in \cite{Zhang09}. Clearly, Lemma 3.4 is a
generalization of their result.

Next, we discuss the error bound for the EVLCP
($\mathbf{A},\mathbf{q}$).

Assume that $x^{\ast}$ is the unique solution of the EVLCP
($\mathbf{A},\mathbf{q}$). Based on Lemma 3.1, we set
\begin{equation}\label{eq:33}
a_{j}=A_{j}x+q_{j}, b_{j}=A_{j}x^{\ast}+q_{j}, j=0,1,\ldots,k.
\end{equation}
Substituting (\ref{eq:33}) into (\ref{eq:31}) yields
\begin{equation}\label{eq:34}
r(x)=\min_{0\leq j\leq k}\{A_{j}x+q_{j}\}-\min_{0\leq j\leq
k}\{A_{j}x^{\ast}+q_{j}\}=(D_{0}A_{0}+D_{1}A_{1}+\ldots+D_{k}A_{k})(x-x^{\ast}),
\end{equation}
where  $D_{j}=\mbox{diag}(d_{j})$ with  $d_{j}\in [0,1]^{n}$
$(j=0,1,\ldots,k)$ are nonnegative diagonal matrices and
$\sum_{j=0}^{k}D_{j}=I$.

Combining (\ref{eq:34}) with Lemma 3.4, we immediately gain the
upper global error bound of the EVLCP ($\mathbf{A},\mathbf{q}$)
under the row $\mathcal{W}$-property.

\begin{theorem} Let $D_{0}, D_{1},\ldots, D_{k}$ satisfy the
conditions of Lemma 3.4. If $\mathbf{A}=(A_{0},A_{1},\ldots, A_{k})$
has the row $\mathcal{W}$-property, then for any $x\in
\mathbb{R}^{n}$,
\begin{equation}\label{eq:35}
\|x-x^{\ast}\|\leq\max\|(D_{0}A_{0}+D_{1}A_{1}+\ldots+D_{k}A_{k})^{-1}\|\|r(x)\|.
\end{equation}
\end{theorem}

\textbf{Remark 3.2} When $A_{0}=I$ and $k=1$ in (\ref{eq:35}),
Theorem 3.1 reduces to (\ref{eq:13}), see Eq. (2.3) on page 516 in
\cite{Chen06} as well. Further, when $k=1$ in (\ref{eq:35}), Theorem
3.1 reduces to Theorem 3.1 in \cite{Zhang09}. In addition, comparing
(\ref{eq:35}) with (\ref{eq:15}), the former advantage over the
latter is that the former  no longer requires the row rearrangement
of the matrix $\mathbf{A}$. That is to say, (\ref{eq:35})
successfully avoids the row rearrangement of the matrix
$\mathbf{A}$. From the view of the calculation time, in general, Eq.
(\ref{eq:15}) requires $(\frac{(k+1)k}{2})^{n}$ times, Eq.
(\ref{eq:35}) requires only 1 time. Compared with the error bound
(\ref{eq:15}) by the row rearrangement technique, our new error
bound (\ref{eq:35}) greatly reduces the computation workload in a
way.

In addition, we also obtain the lower global error bound for the
EVLCP ($\mathbf{A},\mathbf{q}$) under the row
$\mathcal{W}$-property. By making use of (\ref{eq:34}), we easily
find that
\[
\|x-x^{\ast}\|\geq\frac{\|r(x)\|}{\max\|D_{0}A_{0}+D_{1}A_{1}+\ldots+D_{k}A_{k}\|}.
\]

Theorem 3.1 provides a new result for the upper global  error bound
of the EVLCP ($\mathbf{A},\mathbf{q}$) under the row
$\mathcal{W}$-property,  whereas, it contains arbitrary nonnegative
diagonal matrices $D_{0}, D_{1},\ldots, D_{k}$ such that it seems be
an impossible task to find
$\max\|(D_{0}A_{0}+D_{1}A_{1}+\ldots+D_{k}A_{k})^{-1}\|$ under
certain conditions. To turn around this negative situation, an
effective approach is to limit the range of $\mathbf{A}$. That is to
say, once we are able to choose the suitable $\mathbf{A}$, Theorem
3.1 is able to get rid of this unfavorable situation, on condition
that the error can be permitted.

In the following, we consider the upper global  error bound of the
EVLCP ($\mathbf{A},\mathbf{q}$) from two aspects: (I)
$\mathbf{A}=(A_{0},A_{1},\ldots, A_{k})$ with the diagonal part of
every matrix $A_{j}$  being positive; (II)
$\mathbf{A}=(A_{0},A_{1},\ldots, A_{k})$ with every matrix $A_{j}$
being strictly diagonally dominant.
% For case (I), interestingly, we
%obtain the upper global error bound the same as these in
%\cite{Zhang09}, see Theorem 3.2. Not only that, our proof is simpler
%than the proof of Theorem 4.3 in \cite{Zhang09}. For case (II), we
%give a new result, see Theorem 3.3, which is sharper than Theorem
%4.4 in  \cite{Zhang09} in a way.
Under these two cases, interestingly, we obtain the upper global
error bound the same as these in \cite{Zhang09}, see Theorem 3.2 and
Theorem 3.3. Not only that, The proof of Theorem 3.2 is simpler than
the proof of Theorem 4.3 in \cite{Zhang09}.

\begin{theorem} Let  $A_{j}=\wedge_{j}-C_{j}$ in $\mathbf{A}$ with
$\wedge_{j}>0$, where $\wedge_{j}$ is the diagonal part of $A_{j}$,
j=0,1,\ldots,k. If $A_{j}=\wedge_{j}-C_{j}$ satisfy
\begin{equation}\label{eq:36}
\rho(\max_{0\leq j\leq k}\{\wedge^{-1}_{j}|C_{j}|\})<1,
\end{equation}
then $\mathbf{A}=(A_{0},A_{1},\ldots, A_{k})$ has the row
$\mathcal{W}$-property and
\begin{equation}\label{eq:37}
\max\|(D_{0}A_{0}+D_{1}A_{1}+\ldots+D_{k}A_{k})^{-1}\|\leq\|[I-\max_{0\leq
j\leq k}\{\wedge^{-1}_{j}|C_{j}|\}]^{-1}\max_{0\leq j\leq
k}\{\wedge^{-1}_{j}\}\|.
\end{equation}
\end{theorem}
\textbf{Proof.} Let $V=\sum_{j=0}^{k}D_{j}\wedge_{j}$ and
$U=\sum_{j=0}^{k}D_{j}C_{j}$. Then
\[
(D_{0}A_{0}+D_{1}A_{1}+\ldots+D_{k}A_{k})^{-1}=(V-U)^{-1}=(I-V^{-1}U)^{-1}V^{-1}.
\]
Under the condition (\ref{eq:36}), together with Lemma 3.2, we have
\[
V^{-1}\leq\max_{0\leq j\leq k}\{\wedge^{-1}_{j}\}, \ V^{-1}|U|\leq
\max_{0\leq j\leq k}\{\wedge^{-1}_{j}|C_{j}|\},
\]
and
\begin{align*}
|(I-V^{-1}U)^{-1}|=&|I+(V^{-1}U)+(V^{-1}U)^{2}+...|\\
\leq&I+(V^{-1}|U|)+(V^{-1}|U|)^{2}+...\\
\leq&I+(\max_{0\leq j\leq k}\{\wedge^{-1}_{j}|C_{j}|\})+(\max_{0\leq
j\leq k}\{\wedge^{-1}_{j}|C_{j}|\})^{2}+...\\
=&[I-\max_{0\leq j\leq k}\{\wedge^{-1}_{j}|C_{j}|\}]^{-1}.
\end{align*}
So,
\[
|(I-V^{-1}U)^{-1}V^{-1}|\leq[I-\max_{0\leq j\leq
k}\{\wedge^{-1}_{j}|C_{j}|\}]^{-1}\max_{0\leq j\leq
k}\{\wedge^{-1}_{j}\}.
\]
Hence,
\begin{align*}
\|(I-V^{-1}U)^{-1}V^{-1}\|\leq&\||(I-V^{-1}U)^{-1}V^{-1}|\|\\
\leq&\|[I-\max_{0\leq j\leq
k}\{\wedge^{-1}_{j}|C_{j}|\}]^{-1}\max_{0\leq j\leq
k}\{\wedge^{-1}_{j}\}\|,
\end{align*}
which implies that (\ref{eq:37}) is true. By Lemma 3.4, it is easy
to see that $\mathbf{A}=(A_{0},A_{1},\ldots, A_{k})$ has the row
$\mathcal{W}$-property. $\hfill{} \Box$

In fact, the process of the proof of Theorem 3.2 is also suitable
for the proof of Theorem 2.1 in \cite{Chen06} and Theorem 3.2 in
\cite{Zhang09}.

\begin{theorem} Let $\mathbf{A}=(A_{0},A_{1},\ldots,
A_{k})$ with every matrix $A_{j}$  being strictly row diagonally
dominant and sign($\wedge_{0}$)=sign($\wedge_{j}$), where
$\wedge_{j}$ is a diagonal part of $A_{j}$, $j=0,1,\ldots,k$. Then
$\mathbf{A}=(A_{0},A_{1},\ldots, A_{k})$ has the row
$\mathcal{W}$-property and
\begin{equation}\label{eq:38}
\max\|(D_{0}A_{0}+D_{1}A_{1}+\ldots+D_{k}A_{k})^{-1}\|_{\infty}\leq\frac{1}{\min_{i\in
N}\min((\langle A_{0}\rangle e)_{i},(\langle A_{1}\rangle
e)_{i},\ldots, (\langle A_{k}\rangle e)_{i})}.
\end{equation}
%\begin{equation}\label{eq:39}
%\max\|(D_{0}A_{0}+D_{1}A_{1}+\ldots+D_{k}A_{k})^{-1}\|_{\infty}\leq\|(\min_{0\leq
%j\leq k}(\langle A_{j}\rangle))^{-1}\|_{\infty}.
%\end{equation}
\end{theorem}

Here, the proof of Theorem 3.3 is omitted, one can see the proof of
Theorem 3.3 in \cite{Zhang09} for more details.

Comparing Theorem 4.4 in \cite{Zhang09} with Theorem 3.3, the latter
no longer requires ``each $i\in N$, $(A_{j})_{ii}(A_{l})_{ii}>0$,
for any $j<l\in\{0,1,\ldots,k\}$'', just needs to keep $\wedge_{j}$
the same sign.

In \cite{Xiu02}, Xiu and Zhang extended Eq. (\ref{eq:12}) to the
EVLCP ($\mathbf{A},\mathbf{q}$) for $k=1$ under the row
$\mathcal{W}$-property:
\begin{equation}\label{eq:39}
\|x-x^{\ast}\|_{\infty}\leq\frac{\|A_{0}+A_{1}\|_{\infty}}{\alpha\{A_{0},A_{1}\}}\|r(x)\|_{\infty},
\ \mbox{for any }\ x\in \mathbb{R}^{n},
\end{equation}
where
\[
\alpha\{A_{0},A_{1}\}:=\min_{\|x\|_{\infty}=1}\big\{\max_{1\leq
i\leq n}(A_{0}x)_{i}(A_{1}x)_{i}\big\}.
\]
In the sequel, what we're interested in is whether we will extend
(\ref{eq:39}) to the EVLCP ($\mathbf{A},\mathbf{q}$) under the row
$\mathcal{W}$-property, meanwhile, and  avoid the row rearrangement
technique as well. To this end, we require the following lemma, see
Lemma 3.5.

\begin{lemma}
If matrix
$\mathbf{\bar{A}}=(A_{0},D_{1}A_{1}+D_{2}A_{2}+\ldots+D_{k}A_{k})$
has the row $\mathcal{W}$-property, where matrices
$D_{j}=\mbox{diag}(d_{j})$ with  $d_{j}\in [0,1]^{n}$
$(j=1,\ldots,k)$ are arbitrary  nonnegative diagonal matrices and
$\sum_{j=1}^{k}D_{j}=I$, if and only if matrix
$\mathbf{A}=(A_{0},A_{1},\ldots, A_{k})$ has the row
$\mathcal{W}$-property.
\end{lemma}
\textbf{Proof.} ($\Rightarrow$) If
$\mathbf{\bar{A}}=(A_{0},D_{1}A_{1}+D_{2}A_{2}+\ldots+D_{k}A_{k})$
has the row $\mathcal{W}$-property,  matrices
$D_{j}=\mbox{diag}(d_{j})$ with  $d_{j}\in [0,1]^{n}$
$(j=1,\ldots,k)$ are arbitrary  nonnegative diagonal matrices and
$\sum_{j=1}^{k}D_{j}=I$, then, by using Lemma 3.4, for arbitrary
nonnegative diagonal matrix $\hat{D}=\mbox{diag}(\hat{d})$ with
$\hat{d}\in [0,1]^{n}$, matrix
\[
(I-\hat{D})A_{0}+\hat{D}(D_{1}A_{1}+D_{2}A_{2}+\ldots+D_{k}A_{k})
\]is nonsingular.
For arbitrary nonnegative diagonal matrix
$\bar{D}_{j}=\mbox{diag}(\bar{d}_{j})$ with  $\bar{d}_{j}\in
[0,1]^{n}$ $(j=0,1,\ldots,k)$ and $\sum_{j=0}^{k}\bar{D}_{j}=I$, we
set $\hat{D}=\mbox{diag}(\hat{d})=I-\bar{D}_{0}$, i.e.,
$\hat{d}_{i}=1-(\bar{d}_{0})_{i}$, $\hat{d}\in [0,1]^{n}$. Let
$D_{l}=\mbox{diag}(d_{l})$, $l=1,2,\ldots,k$, where
$(d_{l})_{i}=\frac{(\bar{d}_{l})_{i}}{\hat{d}_{i}}$ if
$\hat{d}_{i}\neq 0$; $(d_{l})_{i}=\frac{1}{k}$ if $\hat{d}_{i}=0$.
Further,
\[
\bar{D}_{0}A_{0}+\bar{D}_{1}A_{1}+\ldots
+\bar{D}_{k}A_{k}=(I-\hat{D})A_{0}+\hat{D}(D_{1}A_{1}+D_{2}A_{2}+\ldots+D_{k}A_{k}),
\] and $\sum_{j=1}^{k}D_{j}=I$. So matrix
\[\bar{D}_{0}A_{0}+\bar{D}_{1}A_{1}+\ldots+\bar{D}_{k}A_{k}\] is
nonsingular. By using Lemma 3.4, $\mathbf{A}$ has the row
$\mathcal{W}$-property.

($\Leftarrow$) If $\mathbf{A}$  has the row $\mathcal{W}$-property,
then for arbitrary nonnegative diagonal matrix
$\bar{D}_{j}=\mbox{diag}(\bar{d}_{j})$ with  $\bar{d}_{j}\in
[0,1]^{n}$ $(j=0,1,\ldots,k)$ and $\sum_{j=0}^{k}\bar{D}_{j}=I$,
matrix
\[
\bar{D}_{0}A_{0}+\bar{D}_{1}A_{1}+\ldots+\bar{D}_{k}A_{k}
\]
is nonsingular. In the following, we  prove that
$\mathbf{\bar{A}}=(A_{0},D_{1}A_{1}+D_{2}A_{2}+\ldots+D_{k}A_{k})$
has the row $\mathcal{W}$-property, where matrices
$D_{j}=\mbox{diag}(d_{j})$ with  $d_{j}\in [0,1]^{n}$
$(j=1,\ldots,k)$ are arbitrary  nonnegative diagonal matrices and
$\sum_{j=1}^{k}D_{j}=I$. In fact, for arbitrary nonnegative diagonal
matrix $\hat{D}=\mbox{diag}(\hat{d})$ with $\hat{d}\in [0,1]^{n}$,
\[
(I-\hat{D})A_{0}+\hat{D}(D_{1}A_{1}+D_{2}A_{2}+\ldots+D_{k}A_{k})=(I-\hat{D})A_{0}+\hat{D}D_{1}A_{1}+\hat{D}D_{2}A_{2}+\ldots+\hat{D}D_{k}A_{k}.
\]
Let $\bar{D}_{0}=I-\hat{D}, \bar{D}_{1}=\hat{D}D_{1},\ldots,$
$\bar{D}_{k}=\hat{D}D_{k}$. Then  $\bar{d}_{j}\in [0,1]^{n}$
$(j=0,1,\ldots,k)$, $\sum_{j=0}^{k}\bar{D}_{j}=I$, and matrix
\[
(I-\hat{D})A_{0}+\hat{D}(D_{1}A_{1}+D_{2}A_{2}+\ldots+D_{k}A_{k})=\bar{D}_{0}A_{0}+\bar{D}_{1}A_{1}+\ldots+\bar{D}_{k}A_{k}
\]
is nonsingular. By using Lemma 3.4 again, $\mathbf{\bar{A}}$  has
the row $\mathcal{W}$-property. $\hfill{} \Box$

From Lemma 3.5, together with Lemma 3.1, for the (\ref{eq:39}) type
error bound, we give a general result for the EVLCP
($\mathbf{A},\mathbf{q}$) under the row $\mathcal{W}$-property.

\begin{theorem} Let matrix $\mathbf{A}=(A_{0},A_{1},\ldots, A_{k})$ have the row
$\mathcal{W}$-property. Then for any $x\in \mathbb{R}^{n}$,
\begin{equation}\label{eq:310}
\|x-x^{\ast}\|_{\infty}\leq\min_{D_{1},\ldots,
D_{k}}\bigg\{\frac{\|A_{0}+\sum_{j=1}^{k}D_{j}A_{j}\|_{\infty}}{\alpha\{A_{0},
A_{1}, \ldots,A_{k}\}}\bigg\}\|r(x)\|_{\infty},
\end{equation}
where matrices $D_{j}=\mbox{diag}(d_{j})$ with  $d_{j}\in [0,1]^{n}$
$(j=1,\ldots,k)$ are arbitrary  nonnegative diagonal matrices and
$\sum_{j=1}^{k}D_{j}=I$, and
\[
\alpha\{A_{0}, A_{1}, \ldots,A_{k}\}
=\min_{\|x\|_{\infty}=1}\Bigg\{\max_{1\leq i\leq
n}(A_{0}x)_{i}\bigg(\bigg(\sum_{j=1}^{k}D_{j}A_{j}\bigg)x\bigg)_{i}\Bigg\}.
\]
\end{theorem}
\textbf{Proof.} Since matrix $\mathbf{A}=(A_{0},A_{1},\ldots,
A_{k})$ has the row $\mathcal{W}$-property, from Lemma 3.5, matrix
$\mathbf{\bar{A}}=(A_{0},D_{1}A_{1}+D_{2}A_{2}+\ldots+D_{k}A_{k})$
has the row $\mathcal{W}$-property, where  matrices
$D_{j}=\mbox{diag}(d_{j})$ with  $d_{j}\in [0,1]^{n}$
$(j=1,\ldots,k)$ are arbitrary  nonnegative diagonal matrices and
$\sum_{j=1}^{k}D_{j}=I$.  As done in \cite{Xiu02}, we introduce a
quantity $\alpha\{A_{0}, A_{1}, \ldots,A_{k}\}$ below
\[
\alpha\{A_{0}, A_{1},
\ldots,A_{k}\}=\min_{\|x\|_{\infty}=1}\Bigg\{\max_{1\leq i\leq
n}(A_{0}x)_{i}\bigg(\bigg(\sum_{j=1}^{k}D_{j}A_{j}\bigg)x\bigg)_{i}\Bigg\}.
\]
Further, we take
\[
s(x) = (A_{0}x + q_{0}) - r(x)\  \mbox{and}\ t(x) =\min_{1\leq j\leq
k}\{A_{j}x + q_{j}\}-r(x).
\]
Then
\[
s(x)\geq0, t(x)\geq0, s(x)^{T}t(x)=0.
\]
Thus, for any $x\in \mathbb{R}^{n}$ and each $i = 1,2,. . , n$,
together with Lemma 3.1, we obtain
\begin{align*}
0\geq &(s_{i}(x)-s_{i}(x^{\ast}))(t_{i}(x)-t_{i}(x^{\ast})) \\
=& ((A(x - x^{*}))_{i} - r_{i}(x) +
r_{i}(x^{\ast}))(((\sum_{j=1}^{k}D_{j}A_{j})(x -
x^{\ast}))_{i}-r_{i}(x) + r_{i}(x^{\ast}))\\
\geq &((A(x - x^{\ast}))_{i}((\sum_{j=1}^{k}D_{j}A_{j})(x -
x^{\ast}))_{i}\\
&- (r_{i}(x) - r_{i}(x^{\ast}))((A +
\sum_{j=1}^{k}D_{j}A_{j})(x - x^{\ast}))_{i}\\
\geq &((A(x - x^{\ast}))_{i}((\sum_{j=1}^{k}D_{j}A_{j})(x -
x^{\ast}))_{i}\\
&- \|r(x)\|_{\infty}\cdot\|A_{0}
+\sum_{j=1}^{k}D_{j}A_{j}\|_{\infty}\cdot\|x - x^{\ast}\|_{\infty},
\end{align*}
from which it follows that
\begin{align*}
\|A_{0}
+\sum_{j=1}^{k}D_{j}A_{j}\|_{\infty}\cdot\|r(x)\|_{\infty}\cdot\|x -
x^{\ast}\|_{\infty}\geq& \max_{1\leq i\leq n}((A(x -
x^{\ast}))_{i}((\sum_{j=1}^{k}D_{j}A_{j})(x -
x^{\ast}))_{i}\\
\geq&\alpha\{A_{0}, A_{1}, \ldots,A_{k}\}\|x - x^{\ast}\|_{\infty}.
\end{align*}
This yields the desired inequality (\ref{eq:310}). $\hfill{} \Box$

\section{Numerical examples}
Since the advantages and disadvantages of Theorem 3.2 and Theorem
3.3 for numerical examples have been presented in \cite{Zhang09}, we
here compare Theorem 3.1 with Eq. (\ref{eq:15}), also see Theorem
4.1 in \cite{Zhang09}. Numerical examples used by us are from two
aspects: (1) on the one hand, we still adopt a example in
\cite{Zhang09}; on the other hand, we list some new examples.

\textbf{Example 4.1} \cite{Zhang09} Let $\mathbf{A}=(A_{0}, A_{1},
A_{2})$, where
\[
A_{0}=\left[\begin{array}{cc}
1&1\\
-1&1
\end{array}\right],\ A_{1}=\left[\begin{array}{cc}
2&1\\
-1&1
\end{array}\right], \ A_{2}=\left[\begin{array}{cc}
1&3\\
-1&1
\end{array}\right].
\]
Since $(A_{0})_{2\cdot}=(A_{1})_{2\cdot}=(A_{2})_{2\cdot}$, this
particularity avoids the row rearrangement. Even so, by  Eq.
(\ref{eq:15}), their calculations take \textbf{3} times:
\begin{align*}
\mu_{1}=&\max_{d\in [0,1]^{n}}\|((I-D)A_{0}+DA_{1})^{-1}\|_{\infty}=1,\\
\mu_{2}=&\max_{d\in [0,1]^{n}}\|((I-D)A_{0}+DA_{2})^{-1}\|_{\infty}=1,\\
\mu_{3}=&\max_{d\in
[0,1]^{n}}\|((I-D)A_{1}+DA_{2})^{-1}\|_{\infty}=1.
\end{align*}
Based on this,
\[
\alpha_{\infty}(\textbf{A})=\max\{\mu_{1},\mu_{2},\mu_{3}\}=1.
\]
Directly using Theorem 3.1, for $d_{1}+d_{2}\leq1$\ \mbox{with}\
$d_{1},d_{2}\in [0,1]$,  we calculate
\[
\max_{d\in
[0,1]^{n}}\|(D_{0}A_{0}+D_{1}A_{1}+D_{2}A_{2})^{-1}\|_{\infty}=\max\{\frac{2+2d_{2}}{2+2d_{2}+d_{1}},
\frac{2+d_{1}}{2+d_{1}+2d_{2}}\}=1=\alpha_{\infty}(\textbf{A}).
\]

Whereas, once the system matrix occurs the row rearrangement, the
calculation times by Eq. (\ref{eq:15}) increases dramatically, see
Example 2.1 and Example 4.3.

\textbf{Example 4.2} Here, Example 2.1 is used as Example 4.2 to
investigate Theorem 3.1. Because the main diagonal elements of
$A_{0}, A_{1}, A_{2}$ are positive, the second element in the first
row is non-negative, and the first element in the second row is
non-positive, so, for matrix $D_{0}A_{0}+D_{1}A_{1}+D_{2}A_{2}$, its
main diagonal elements are positive, the second element in the first
row is non-negative, and the first element in the second row is
non-positive for any $D_{j}=\mbox{diag}(d_{j})$ with  $d_{j}\in
[0,1]^{n}$ $(j=0,1,2)$ being nonnegative diagonal matrices and
$\sum_{j=0}^{2}D_{j}=I$. According to the definition of the
determinant, we can get $\det(D_{0}A_{0}+D_{1}A_{1}+D_{2}A_{2})>0$,
i.e., $D_{0}A_{0 }+D_{1}A_{1}+D_{2}A_{2}$ is nonsingular. From Lemma
3.4, the block matrix $\mathbf{A}=(A_{0}, A_{1}, A_ {2})$ has row
$\mathcal{W}$-property.

Directly using Theorem 3.1, for $D_{j}=\mbox{diag}(d_{j})$ with
$d_{j}\in [0,1]^{n}$ $(j=0,1,2)$ are nonnegative diagonal matrices
and $\sum_{j=0}^{2}D_{j}=I$, we calculate
\[
\max_{d\in
[0,1]^{n}}\|(D_{0}A_{0}+D_{1}A_{1}+D_{2}A_{2})^{-1}\|_{\infty}=\max\left\{w_{1},
w_{2}\right\}=3,
\]
where
\[
w_{1}=\frac{2-(d_{1})_{11}}{1+(d_{2})_{11}+(1-(d_{1})_{11})(1+(d_{1})_{22}-(d_{2})_{22})},
\]
and
\[
w_{2}=\frac{2+(d_{1})_{22}}{1+(d_{2})_{11}+(1-(d_{1})_{11})(1+(d_{1})_{22}-(d_{2})_{22})}.
\]
When using Eq. (\ref{eq:15}), although
$\alpha_{\infty}(\textbf{A})=3$, see Example 2.1, which is equal to
our this result, its calculation times are \textbf{9}.

It is easy to check that $\mathbf{A}_0$, $\mathbf{A}_1$ in Example
4.2 are not strictly diagonal dominant, and that
$$\rho(\max\{\Lambda_0^{-1}|C_0|,\Lambda_1^{-1}|C_1|,\Lambda_2^{-1}|C_2|\})=1.4142>1.$$
Hence  the error bounds in  Theorem 3.2, Theorem 3.3, Theorem 4.3
and Theorem 4.4 of \cite{Zhang09} cannot work for this case.

\textbf{Example 4.3}  Let $\mathbf{A}=(A_{0}, A_{1}, A_{2}, A_{3})$,
where
\[
A_{0}=\left[\begin{array}{cc}
1&1\\
-1&1
\end{array}\right],\ A_{1}=\left[\begin{array}{cc}
1&0\\
-2&1
\end{array}\right], \ A_{2}=\left[\begin{array}{cc}
4&1\\
0&1
\end{array}\right], \ A_{3}=\left[\begin{array}{cc}
2&0\\
-3&1
\end{array}\right]
\]
With the same analysis method as in Example 4.3, it is easy to get
that $\mathbf{A}=(A_{0}, A_{1}, A_{2}, A_{3})$ has row
$\mathcal{W}$-property. Directly using Theorem 3.1, for $D_{0},
D_{1}, D_{2}, D_{3}$ are nonnegative diagonal matrices with
$d_{j}\in [0,1]^{n}$ $(j=0,1,2,3)$ and $\sum_{j=0}^{3}D_{j}=I$, we
calculate
\[
\max_{d\in
[0,1]^{n}}\|(D_{0}A_{0}+D_{1}A_{1}+D_{2}A_{2}+D_{3}A_{3})^{-1}\|_{\infty}=\max\left\{\bar{w}_{1},
\bar{w}_{2}\right\}=4=\alpha_{\infty}(\textbf{A}).
\]
where
\[
\bar{w}_{1}=\frac{2-(d_{1})_{22}-(d_{2})_{22}+2(d_{3})_{22}}{1+3(d_{2})_{11}+(d_{3})_{11}+(1+(d_{1})_{22}-(d_{2})_{22}+2(d_{3})_{22})(1-(d_{1})_{11}-(d_{3})_{11})}
\]
and
\[
\bar{w}_{2}=\frac{2-(d_{1})_{11}+3(d_{2})_{11}}{1+3(d_{2})_{11}+(d_{3})_{11}+(1+(d_{1})_{22}-(d_{2})_{22}+2(d_{3})_{22})(1-(d_{1})_{11}-(d_{3})_{11})}.
\]
Of course, we can adopt Eq. (\ref{eq:15}) to obtain
$\alpha_{\infty}(\textbf{A})=4$, whereas, its calculation times are
\textbf{36}. Here, the computational process of Eq. (\ref{eq:15}) is
omitted.

It is easy to check that $\mathbf{A}_0$, $\mathbf{A}_1$,
$\mathbf{A}_3$ are not strictly diagonal dominant, and that
\[\rho(\max\{\Lambda_0^{-1}|C_0|,\Lambda_1^{-1}|C_1|,\Lambda_2^{-1}|C_2|,
\Lambda_3^{-1}|C_3|\})=2.3028>1.\] For Example 4.3,  the error
bounds in  Theorem 3.2, Theorem 3.3, Theorem 4.3 and Theorem 4.4 of
\cite{Zhang09} cannot work for this case as well.

From Example 2.1 and Example 4.3, it is easy to find that although
the latter has one more matrix than the former, the calculation
times of the latter are more \textbf{27} times than the former.

\section{Conclusion}
In this paper, by introducing a general equivalent form of the
minimum function, some new error bounds for the
EVLCP($\mathbf{A},\mathbf{q}$) under the row $\mathcal{W}$-property
have been presented. These new error bounds not only cover some
existing results, but also keep away from the row rearrangement of
the system matrix. In addition, with respect to the row
$\mathcal{W}$-property, we also obtain two new sufficient and
necessary conditions. Finally, by numerical examples, we show that
the new error bound is feasible. Compared with the error bound by
the row rearrangement technique, our new error bound greatly reduces
the computation workload in a way.

In addition, by a lot of numerical experiments, we find that the
error bound (\ref{eq:35}) is equal to the error bound (\ref{eq:15}).
There exists an interested problem how to prove the equality of
both, which is necessary for us to study it in the future.

%\section*{Acknowledgments}
%
%The authors would like to thank two anonymous referees for providing
%helpful suggestions, which greatly improved the paper.

\end{document}